\input amstex 
\input xy
\xyoption{all}
 
 \documentstyle{amsppt}

\document
\magnification=1200
\NoBlackBoxes
\nologo
\vsize18cm
\centerline{\bf  Set-theoretical solutions to the Yang-Baxter Relation }
 \medskip
\centerline{\bf from factorization of matrix polynomials and 
$\theta$-functions}   
\bigskip
\centerline{\bf Alexander Odesskii}
\bigskip
\centerline{\bf Introduction}

\medskip

The Yang-Baxter relation plays a central role in two-dimensional Quantum Field Theory. 
This relation involves a linear operator $R : V\otimes V\to V\otimes V$, 
where $V$ is a vector space, and has the form

$$ R^{12}R^{13}R^{23}=R^{23}R^{13}R^{12}$$ 
in $End(V\otimes V\otimes V)$, where $R^{ij}$ means $R$ acting in the
$i$-th and $j$-th components. In the paper [12] V. Drinfeld suggested to study
set-theoretical solutions of 
this relation, i.e. solutions given by a map $R : X\times X\to X\times X$,
where $X$ is a given set. 
Moreover, if $X$ is an algebraic manifold, then $R$ may be a rational
map. The general theory of 
set-theoretical solutions to the quantum Yang-Baxter relation was developed
in 
[11, 13].
Various examples were constructed in [10, 11, 13]. In this paper we
construct such solutions 
from decompositions of matrix polynomials and $\theta$-functions. These
solutions arise 
from the decompositions "in different order". We also construct a "local
action 
of the symmetric group" in these cases, generalizations of the action of
the symmetric group $S N$ on 
$X^N$ given by the set-theoretical solution. The structure of the paper is 
as follows. In {\bf \S1} we give basic definitions. In {\bf \S2} we
introduce a set-theoretical solution arising from the factorization of 
matrix polynomials. In {\bf \S3} we introduce a set-theoretical solution 
arising from matrix $\theta$-functions.

For a given set-theoretical solution of the quantum Yang-Baxter relation
one can define a twisted Yang-Baxter relation with the set of spectral 
parameters $X$ (see [14] and (3) of this paper). The corresponding twisted 
$R$-matrix describes a scattering of two "particles" such that the spectral 
parameters change after scattering according to a given set-theoretical 
solution. Moreover, one can define a generalized star-triangle relation for 
a given local action of the symmetric group (see [14]). The examples of 
twisted $R$-matrices as well as the solutions of the generalized
star-triangle relation were found in [4] as intertwiners of cyclic 
representations and their tensor products of the algebra of monodromy 
matrices of the six-vertex model at roots of unity [3]. These solutions 
are natural generalizations of the one from the chiral Potts model
[1,2,3]. Other examples were found in [9] for the relativistic Toda chain. 
One can obtain various solutions of the twisted Yang-Baxter and
star-triangle relations by calculating the intertwiners of the 
representations of the algebras of monodromy matrices at roots of unity for 
other trigonometric and elliptic $R$-matrices.

\newpage

\centerline{\bf \S1.  Basic definitions}
\medskip

Let $U$ be a complex manifold, $\mu: U\times U\to U\times U$ be a
birational automorphism of $U\times U$. We will use a notation: 
$\mu(u,v)=(\varphi(u,v),\psi(u,v))$ where $u, v \in U$.
Here $\varphi$ and $\psi$ are meromorphic functions from $U\times U$ to $U$.

Let us introduce the following birational automorphisms of $U\times U\times
U$: $\sigma_1=\mu\times \text{id}$ and $\sigma_2=\text{id}\times\mu$. 
We have: $\sigma_1(u,v,w)=(\varphi(u,v),\psi(u,v),w)$ and $\sigma_2(u,v,w)=
(u,\varphi(v,w),\psi(v,w))$.

{\bf Definition} {\it We call a map $\mu$ a twisted transposition if 
the automorphisms $\sigma_1$ and $\sigma_2$ satisfy the following relations:}

$$\sigma_1^2=\sigma_2^2=\text{id},
\text{     } \sigma_1\sigma_2\sigma_1=\sigma_2\sigma_1\sigma_2 \eqno(1)$$

If $\mu$ is a twisted transposition, then for each $N\in \Bbb N$ we 
have a birational action of the symmetric group $S_N$ on the manifold 
$U^N$ such that the transposition $(i,i+1)$ acts by an automorphism 
$\sigma_i=\text{id}^{i-1}\times\mu\times \text{id}^{N-i-1}$. So we have 
$\sigma_i(u_1,\dots,u_N)=(u_1,\dots,\varphi(u_i,u_{i+1}),\psi(u_i,u_{i+1}),
\dots,u_N)$. It is clear that the relations (1) are equivalent to the
defining relations in the group $S_N$: $\sigma_i^2=e, \sigma_i
\sigma_{i+1}\sigma_i=\sigma_{i+1}\sigma_i\sigma_{i+1}, \sigma_i\sigma_j=
\sigma_j\sigma_i$ for $|i-j|>1$.

It is easy to check that the relations (1) are equivalent to the following 
functional equations for $\varphi$ and $\psi$:

$$\varphi(\varphi(u,v),\psi(u,v))=u,\text{   } \psi(\varphi(u,v),\psi(u,v))=v$$

$$\varphi(u,\varphi(v,w))=\varphi(\varphi(u,v),\varphi(\psi(u,v),w))$$

$$\varphi(\psi(u,\varphi(v,w)),\psi(v,w))=\psi(\varphi(u,v),\varphi(\psi(u,v),w)) \eqno(2)$$

$$\psi(\psi(u,v),w)=\psi(\psi(u,\varphi(v,w)),\psi(v,w))$$

{\bf Remarks 1.} From (2) it follows that for each $N$ the functions $\varphi(u_1,\varphi(u_2,\dots,\varphi(u_N,w)\dots)$ and $\psi(\dots(\psi(w,u_1),u_2)\dots,u_N)$ are invariant with respect to the action of the group $S_N$ on the variables $u_1,\dots,u_N$.

{\bf 2.} Let $\sigma : U\times U\to U\times U$ be the map given by $\sigma(u,v)=(v,u)$. Then $\sigma\mu$ is a set-theoretical solution to the quantum Yang-Baxter relation.

{\bf 3.} Informally one can consider $\sigma\mu$ as an infinite dimensional $R$-matrix in the space of functions. Namely, if we consider the space of meromorphic functions $\{f: U\times U\to \Bbb C\}$ as an "extended tensor square" of the space of meromorphic functions $\{f: U\to\Bbb C\}$, then the linear operator $R_{\sigma\mu}: f\to f\sigma\mu$ (that is $R_{\sigma\mu}f(u,v)=f(\sigma(\mu(u,v)))$) satisfies the usual Yang-Baxter relation.

{\bf Examples 1.} Let $q, q^{-1}: U\to U$ be birational automorphisms such that $qq^{-1}=q^{-1}q=\text{id}$. Then 
$$\mu(u,v)=(q(v),q^{-1}(u))$$
is a twisted transposition.

{\bf 2.} Let $U=\Bbb C$, then the following formula gives a twisted transposition:

$$\mu(u,v)=(1-u+uv, \frac{uv}{1-u+uv})$$

{\bf 3.} Let $U$ be a finite dimensional associative algebra with a unity 
$1\in U$, for example $U=\text{Mat}_m$. Then the following formula gives a twisted transposition:

$$\mu(u,v)=(1-u+uv, (1-u+uv)^{-1}uv)$$

\medskip

Let $V$ be a $n$-dimensional vector space. For each $u\in U$ we denote by $V(u)$ a vector space canonically
isomorphic to $V$. Let $R$ be a meromorphic function from $U\times U$ to $End(V\otimes V)$. We will consider $R(u,v)$ as a linear operator 
$$R(u,v): V(u)\otimes V(v)\to V(\varphi(u,v))\otimes V(\psi(u,v))$$

{\bf Definition} {\it We call $R$ a twisted $R$-matrix (with respect to the twisted transposition $\mu$) if it satisfies the following properties:

1. The composition

$$V(u)\otimes V(v)\to V(\varphi(u,v))\otimes V(\psi(u,v))\to V(u)\otimes V(v)$$
is equal to the identity, that is $R(\varphi(u,v),\psi(u,v))R(u,v)=1$.

2. The following diagram is commutative:
$$\xymatrix{{V(\varphi(u,v))\otimes V(\psi(u,v))\otimes V(w)}\ar[r]^(.42){1\otimes R(\psi(u,v),w)}&{V(\varphi(u,v))\otimes V(\varphi(\psi(u,v)),w))\otimes V(\psi(\psi(u,v),w))}\ar[d]^{R(\varphi(u,v),\varphi(\psi(u,v),w))\otimes 1}\\{V(u)\otimes V(v)\otimes V(w)}\ar[u]_{R(u,v)\otimes 1}\ar[d]^{1\otimes R(v,w)}&{\widetilde V}\\{V(u)\otimes V(\varphi(v,w))\otimes V(\psi(v,w))}\ar[r]^(.42){R(u,\varphi(v,w))\otimes 1}&{V(\varphi(u,\varphi(v,w)))\otimes V(\psi(u,\varphi(v,w)))\otimes V(\psi(v,w))}\ar[u]_{1\otimes R(\psi(u,\varphi(v,w)),\psi(v,w))}}$$

Here $\widetilde V=V(\varphi(u,\varphi(v,w)))\otimes V(\psi(\varphi(u,v),\varphi(\psi(u,v),w)))\otimes V(\psi(\psi(u,v),w))$.

In other words, 

$$R^{12}(\varphi(u,v),\varphi(\psi(u,v),w)))R^{23}(\psi(u,v),w)R^{12}(u,v)=$$
$$R^{23}(\psi(u,\varphi(v,w)),\psi(v,w))R^{12}(u,\varphi(v,w))R^{23}(v,w) \eqno(3)$$

Here $R^{12}=R\otimes 1$ and $R^{23}=1\otimes R$ are linear operators in $V\otimes V\otimes V$.

We call (3) a twisted Yang-Baxter relation.}

Let $\{x_i, i=1\dots,n\}$ be a basis of the linear space $V$, $\{x_i(u)\}$ be the corresponding basis of the linear space $V(u)$. It is clear that the following two linear operators are twisted $R$-matrices for each $\mu$:
$$x_i(u)\otimes x_j(v)\to x_i(\varphi(u,v))\otimes x_j(\psi(u,v))$$ 
$$x_i(u)\otimes x_j(v)\to x_j(\varphi(u,v))\otimes x_i(\psi(u,v))$$

\newpage
\centerline{\bf \S2.  Set-theoretical solution from}
\centerline{\bf factorization of matrix polynomials}
\medskip

For the general theory of matrix polynomials and factorizations see
[7]. For our 
purposes we state results, which may be well known to the experts.

We denote by $S(a)$ the set of eigenvalues of a matrix $a\in Mat_m$. More
generally, 
we denote by $S(a_1,\dots,a_d)$, $a_1,\dots,a_d\in Mat_m$, the set of roots of a polynomial $f(t)=\det(t^d-a_1t^{d-1}+\dots+(-1)^da_d)$. We will consider polynomials with generic coefficients only, so $\#S(a_1,\dots,a_d)=md$.

{\bf Proposition 1.} {\it Let 
$$t^d-a_1t^{d-1}+\dots+(-1)^da_d=(t-b_1)\dots(t-b_d) \eqno(4)$$
for generic matrices $a_1,\dots,a_d\in Mat_m$, then $S(b_i)\cap
S(b_j)=\emptyset$ for 
$i\ne j$ and $S(b_1)\cup\dots\cup S(b_d)=S(a_1,\dots,a_d)$. For each 
decomposition $S(a_1,\dots,a_d)=A_1\cup\dots\cup A_d$, such that 
$\#A_i=m$, $A_i\cap A_j=\emptyset$ $(i\ne j)$ there exists a unique factorization (4) with $S(b_i)=A_i$.}

{\bf Proof} The first statement follows from the equation $\det(t^d-a_1t^{d-1}+\dots+(-1)^da_d)=\det(t-b_1)\dots \det(t-b_d)$.

On the other hand, if we know eigenvalues of $b_1,\dots,b_d$ then we can 
calculate eigenvectors of them. For $\lambda\in S(b_d)$ the corresponding 
eigenvector is a vector $v_{\lambda}$, such that
$({\lambda}^d-a_1{\lambda}^{d-1}+
\dots+(-1)^da_d)v_{\lambda}=0$. If we know all eigenvectors of $b_d$, then 
we can calculate all eigenvectors of $b_{d-1}$ similarly and so on. This 
implies the uniqueness. By our construction of $b_1,\dots,b_d$ the 
determinants of the matrix polynomials in the right hand side and the 
left hand side of (4) have the same sets of roots. Moreover, for each root $\lambda$ the operators represented by these matrix polynomials have the same kernel if we set $t=\lambda$. It implies that these polynomials are equal.

{\bf Proposition 2.} {\it Let $a_1, a_2\in Mat_m$ be generic matrices. Then there exists a unique pair of matrices $b_1, b_2\in Mat_m$ such that $(t-a_1)(t-a_2)=(t-b_1)(t-b_2)$ and $S(b_1)=S(a_2), S(b_2)=S(a_1)$. We have $b_1=a_1+{\Lambda}^{-1}, b_2=a_2-{\Lambda}^{-1}$ where $a_2\Lambda-\Lambda a_1=1$.}

{\bf Proof} If $S(b_2)\cap S(a_2)\ne\emptyset$, then $\det(a_2-b_2)=0$, because $a_2$ and $b_2$ have a common eigenvector. Otherwise, we can put $\Lambda=(a_2-b_2)^{-1}$.

Let $U=Mat_m$. From the propositions 1 and 2 it follows that the formula $\mu(a_1,a_2)=(b_1,b_2)$ gives a twisted transposition, where $b_1+b_2=a_1+a_2$, $b_1b_2=a_1a_2, S(b_1)=S(a_2), S(b_2)=S(a_1)$. We have $\mu(a_1,a_2)=(a_1+{\Lambda}^{-1},a_2-{\Lambda}^{-1})$, where $\Lambda$ is the solution of the linear matrix equation  $a_2\Lambda-\Lambda a_1=1$.

Let $\overline U=\overline{Mat_m}$ be the set of $m\times m$ matrices with different eigenvalues and fixed order of eigenvalues. The proposition 1 gives an action of the symmetric group $S_{mN}$ on the space $\overline U^N$ by birational automorphisms. By definition, for $\sigma\in S_{mN}$, $b_1,\dots,b_N\in \overline U$ we have $\sigma(b_1,\dots,b_N)=(b_1^\prime,\dots,b_N^\prime)$ where $(t-b_1)\dots(t-b_N)=(t-b_1^\prime)\dots(t-b_N^\prime)$ and $\overline S(b_i^\prime)=\sigma\overline S(b_i)$,  $\overline S$ stands for the ordered set of eigenvalues.

This action is local in the following sense. The transposition $(i,i+1)$ for $\alpha m<i<(\alpha+1)m$ acts only inside the $\alpha+1$-th factor of $\bar U^N$ and the transposition $(\alpha m,\alpha m+1)$ acts only inside the product of the $\alpha$-th and the $\alpha+1$-th factors. We have also the twisted transposition $\mu: \overline U\times \overline U\to \overline U\times \overline U$ in this case which is the action of the element $(1,m+1)(2,m+2)\dots(m-1,2m-1)\in S_{2m}$.

{\bf Remark} Let $U$ be the set of matrix polynomials of the form $at+b$, 
where $a,b\in Mat_m, a=(a_{ij}), b=(b_{ij})$ and $a_{ij}=0$ for $i<j$, 
$b_{ij}=0$ for $i>j$. It is possible to define a twisted transposition 
$\mu$ such that for $\mu(f(t),g(t))=(f_1(t),g_1(t))$ we have 
$f(t)g(t)=f_1(t)g_1(t)$, $\det f(t)$ and $\det g_1(t)$ have the same sets of 
roots and the first coefficients of $f(t)$ and $g_1(t)$ have the same
diagonal elements. In [4] we found the solutions of the corresponding 
twisted Yang-Baxter relation (for $m=2$), which is a generalization of the $R$-matrix from chiral Potts model.

\newpage
\centerline{\bf \S3.  Set-theoretical solution from }
\centerline{\bf factorization of matrix $\theta$-functions}
\medskip

Let $\Gamma\subset\Bbb C$ be a lattice generated by 1 and $\tau$ where 
$\text{Im}\tau>0$. We have $\Gamma=\{\alpha+\beta\tau; \alpha,\beta\in\Bbb
Z\}$. 
Let $\varepsilon\in\Bbb C$ be a primitive root of unity of degree $m$. Let $\gamma_1,\gamma_2\in Mat_m$ be $m\times m$ matrices such that $\gamma_1^m=\gamma_2^m=1, \gamma_2\gamma_1=\varepsilon\gamma_1\gamma_2$. We have $\gamma_1v_{\alpha}=\varepsilon^{\alpha}v_{\alpha}, \gamma_2v_{\alpha}=v_{\alpha+1}$ in some basis $\{v_{\alpha}; \alpha\in\Bbb Z/m\Bbb Z\}$ of $\Bbb C^m$. Let us assume that $\{v_1,\dots,v_m\}$ is the standard basis of $\Bbb C^m$. 

We denote by $M\Theta_{n,m,c}(\Gamma)$ for $n,m\in\Bbb N, c\in\Bbb C$ the space of everywhere holomorphic functions $f: \Bbb C\to Mat_m$, which satisfy the following equations:

$$f(z+\frac{1}{m})=\gamma_1^{-1}f(z)\gamma_1$$
$$f(z+\frac{1}{m}\tau)=e^{-2\pi i(mnz-c)}\gamma_2^{-1}f(z)\gamma_2 \eqno(5)$$

{\bf Proposition 3.} {\it $\dim M\Theta_{n,m,c}(\Gamma)=m^2n$ and for each element $f\in M\Theta_{n,m,c}(\Gamma)$ the equation $\det f(z)=0$ has exactly $mn$ zeros modulo $\frac{1}{m}\Gamma$. The sum of these zeros is equal to $mc+\frac{mn}{2}$ modulo $\Gamma$.}

{\bf Proof} For $m=1$ we have the usual $\theta$-functions 
$\Theta_{n,c}(\Gamma)=M\Theta_{n,1,c}(\Gamma)$ and all these 
statements are well known in this case ([8]). One has a basis $\{\theta_{\alpha}(z); \alpha\in\Bbb Z/n\Bbb Z\}$ in the space  $\Theta_{n,c}(\Gamma)$ such that $\theta_{\alpha}(z+\frac{1}{n})=e^{2\pi i\frac{\alpha}{n}}\theta_{\alpha}(z), \theta_{\alpha}(z+\frac{1}{n}\tau)=e^{-2\pi i(z-\frac{n-1}{2n}\tau-\frac{1}{n}c)}\theta_{\alpha+1}(z)$ ([8]). From (5) it follows that $f(z+1)=f(z)$ and $f(z+\tau)=e^{-2\pi i(m^2nz-c_1)}f(z)$ for some $c_1\in\Bbb C$. So the matrix elements of $f(z)$ are $\theta$-functions from the space $\Theta_{m^2n,c_1}(\Gamma)$. We have decomposition $f(z)=\sum_{\alpha}\varphi_{\alpha}\theta_{\alpha}(z)$, where $\varphi_{\alpha}\in Mat_m$ are constant matrices, $\{\theta_{\alpha}\}$ is a basis in the space $\Theta_{m^2n,c_1}(\Gamma)$. Substituting this decomposition in (6) one can calculate the dimension of the space $M\Theta_{n,m,c}(\Gamma)$. We have also $\det f(z+\frac{1}{m})=\det f(z)$ and $\det f(z+\frac{1}{m}\tau)=e^{-2\pi i(m^2nz-mc)}\det f(z)$. From this follows the statement about zeros of the equation $\det f(z)=0$.

{\bf Proposition 4.} {\it For generic complex numbers $\lambda_1,\dots,
\lambda_{mn}$ such that $\lambda_1+\dots+\lambda_{mn}\equiv
mc+\frac{mn}{2}$ 
$\text{mod}\frac{1}{m}\Gamma$ and nonzero vectors $v_1,\dots,v_{mn}\in\Bbb C^m$ there exists a unique up to proportionality element $f(z)\in M\Theta_{n,m,c}(\Gamma)$ such that $\det f(\lambda_{\alpha})=0, f(\lambda_{\alpha})v_{\alpha}=0$ for $1\leqslant\alpha\leqslant mn$.}

{\bf Proof} Considering the decomposition $f(z)=\sum_{\alpha}\varphi_{\alpha}\theta_{\alpha}(z)$, one has the system of linear equations $\{\sum_{\alpha}\theta_{\alpha}(\lambda_{\beta})\varphi_{\alpha}v_{\beta}=0; \beta=1,\dots,mn\}$ for matrix elements of $\{\varphi_{\alpha}\}$. One can see that this system defines $\{\varphi_{\alpha}\}$ uniquely up to proportionality for generic $\lambda_1,\dots,\lambda_{mn}, v_1,\dots,v_{mn}$.

We denote by $S(f)$ the set of zeros of the equation $\det f(z)=0$ modulo $\frac{1}{m}\Gamma$.

{\bf Proposition 5.} {\it Assume that $f(z)\in M\Theta_{n,m,c}(\Gamma)$ is 
a generic element and we have a factorization $f(z)=f_1(z)\dots f_n(z)$, where $f_{\alpha}(z)\in M\Theta_{1,m,c_{\alpha}}(\Gamma)$, $c_1+\dots+c_n=c$. Then $S(f_{\alpha})\cap S(f_{\beta})=\emptyset$ for $\alpha\ne\beta$ and $S(f)=S(f_1)\cup\dots\cup S(f_n)$. For each decomposition $S(f)=A_1\cup\dots\cup A_n$ such that $A_{\alpha}\cap A_{\beta}=\emptyset$ for $\alpha\ne\beta$ and $\#A_{\alpha}=m$ there exists a unique factorization $f(z)=f_1(z)\dots f_n(z)$ up to proportionality of $f_{\alpha}$ such that $S(f_{\alpha})=A_{\alpha}$ for $\alpha=1,\dots,n$.}

{\bf Proof} is similar to the proof of the proposition 1, we just change polynomials by $\theta$-functions.

Let $U_c$ be the projectivisation of the linear space  $M\Theta_{1,m,c}(\Gamma)$ and $U=\bigcup_{c\in\Bbb C}U_c$. We have the following twisted transposition $\mu: U\times U\to U\times U$. By definition $\mu(f,g)=(f_1,g_1)$, where $f(z)g(z)=f_1(z)g_1(z)$ and $S(f_1)=S(g), S(g_1)=S(f)$. 

Let $\overline U$ be the set of elements $f$ from $U$ with a fixed order on $S(f)$. For $f\in\overline U$ let $\overline S(f)$ be the set $S(f)$ with corresponding order. We have a local action of the symmetric group $S_{mN}$ on the space $\overline U^N$. By definition, for $\sigma\in S_{mN}$ we have $\sigma(f_1,\dots,f_N)=(f_1^{\sigma},\dots,f_N^{\sigma})$ where $f_1(z)\dots f_N(z)=f_1^{\sigma}(z)\dots f_N^{\sigma}(z)$ and $\overline S(f_{\alpha}^{\sigma})=\sigma\overline S(f_{\alpha})$ for $1\leqslant\alpha\leqslant N$.

{\bf Remark} It is possible to construct twisted $R$-matrices for this 
twisted transposition $\mu$ as intertwiners of tensor products of 
cyclic representations of the algebra of monodromy matrices for the elliptic Belavin $R$-matrix [5] at the point of finite order (see also [6]). It will be the subject of another paper.

{\bf Acknowledgments}

I am grateful to V.Bazhanov and A.Belavin for useful discussions.

I am grateful to Max-Planck-Institut fur Mathematik, Bonn, where this paper was written, for invitation and very stimulating working atmosphere.

The work is supported partially by RFBR 99-01-01169, RFBR 00-15-96579, CRDF RP1-2254 and INTAS-00-00055.

\newpage
\centerline{\bf References}
\medskip

1. R.J. Baxter, J.H.H. Perk and H. Au-Yang, New solutions of the star-triangle relations for the chiral Potts model, Phys. Lett. A, 128(1988)138-142.

2. V.V. Bazhanov and Yu.G. Stroganov, Chiral Potts model as a descendant of the six vertex model, J. Stat. Phys. 51(1990)799-817.

3. V.O. Tarasov, Cyclic monodromy matrices for the $R$-matrix of the six-vertex model and the chiral Potts model with fixed spin boundary conditions, International Jornal of Modern Physics A, Vol.7, Suppl. 1B(1992)963-975.

4. A.A. Belavin, A.V. Odesskii and R.A. Usmanov, New relations in the algebra of the Baxter $Q$-operators. hep-th/0110126.

5. A.A. Belavin, Discrete groups and integrability of quantum systems, Func. Anal. Appl. 14(1980)18-26.

6. B.L. Feigin and A.V. Odesskii, Sklyanin Elliptic Algebras. The case of point of finite order, Func. Anal. Appl. 29(1995).

7. I. Gelfand,  V. Retakh, R. Wilson, Quadratic-linear algebras associated with decompositions of noncommutative polynomials and noncommutative differential polynomials, q-alg/0002238, to appear in Selecta Mathematica.

8. D. Mumford, Tata lectures on theta. I. With the assistance of C. Musili, M. Nory and M. Stillman. Progress in Mathematics, 28. BirkhDuser Boston, Inc., Boston, MA, 1983. 

9. S. Pakuliak and S. Sergeev, Relativistic Toda chain at root of unity II. Modified $Q$-operator, nlin.SI/0107062.  

10. V. M. Bukhshtaber, The Yang-Baxter transformation, Russian Math. Surveys 53:6 (1998) 1343-1379.

11. P. Etingof, T. Schedler, A. Soloviev, Set-theoretical solutions to the quantum Yang-Baxter equation, q-alg/9801047, Duke Math. J. 1999.

12. V. Drinfeld, Some unsolved problems in quantum group theory, Lecture Notes in Math, 1510, p.1-8.

13. J-H. Lu, M. Yan, Y-C. Zhu, On set-theoretical Yang-Baxter equation, Duke Math. J. 1999.

14. A. Odesskii, Local action of the symmetric group and the twisted Yang-Baxter relation, q-alg/0110268.

\enddocument